\documentclass[10pt]{article}
\usepackage{amsmath,amsthm,amsfonts,amssymb,epsfig,enumerate}
\usepackage{flafter}
\usepackage{algorithm}
\usepackage{algorithmic}

\newtheorem{theorem}{Theorem}[section]

 \newtheorem{lemma}{Lemma}[section]
 \newtheorem{remark}{Remark}
\newtheorem{defn}[theorem]{Definition}
\newtheorem{prop}{Proposition}
\theoremstyle{remark}
\numberwithin{equation}{section}

\newcommand{\norm}[1]{\left\Vert#1\right\Vert}

\usepackage{graphicx}
\usepackage{float}
\usepackage{epstopdf}
\usepackage{multirow}
\usepackage{color}
\begin{document}
\title{On the well-posedness of the time-fractional diffusion equation with Robin boundary condition }
\author{Mengmeng Zhang $^\dag$ \quad  Jijun Liu $^{\dag\; \ddag}\;$\thanks{Corresponding author: Prof. J.J.Liu, email: jjliu@seu.edu.cn}\\
$^\dag$School of Mathematics, Southeast University\\
$^\ddag$Nanjing Center for Applied Mathematics\\
Nanjing, 210096, P.R.China\\
}

\maketitle


\begin{abstract}
The diffusion system with time-fractional order derivative is of great importance mathematically due to the nonlocal property of the fractional order derivative, which can be applied to model the physical phenomena with memory effects. We consider an initial-boundary value problem for the time-fractional diffusion equation with inhomogenous Robin boundary condition. Firstly, we show the unique existence of the weak/strong solution based on the eigenfunction expansions, which ensures the well-posedness of the direct problem. Then, we establish the Hopf lemma for time-fractional diffusion operator, generalizing the counterpart for the classical parabolic equation. Based on this new Hopf lemma, the maximum principles for this time-fractional diffusion are finally proven, which play essential roles for further studying  the uniqueness of the inverse problems corresponding to this system.
\end{abstract}



\section{Introduction}\label{a}

For a bounded domain $\Omega\subset \mathbb{R}^d $ of smooth boundary and $\alpha\in (0,1)$, consider the anomalous diffusion process governed by the differentiation system
\begin{eqnarray}\label{eq1}
\begin{cases}
\partial_{0+}^\alpha u(x,t)-\Delta u(x,t)=f(x)g(t)=:F(x,t), &(x,t)\in \Omega_T:=\Omega\times(0,T],\\
\frac{\partial u(x,t)}{\partial \nu}+\lambda(x)u(x,t)=b(x,t), &(x,t)\in \partial\Omega\times (0,T],\\
u(x,0)=u_0(x), &x\in\overline{\Omega}
\end{cases}
\end{eqnarray}
with time-fractional order derivative,
where $\nu(x)$ is the outward unit normal direction on $\partial \Omega$, and $F(x,t), b(x,t), u_0(x)$ are internal source, boundary source and the initial status of the diffusion process, respectively. The Robin boundary condition in (\ref{eq1}) with impedance coefficient  $\lambda(x)>0$ describes the convection between the solute in a body and one in the ambient environment \cite{WangJG}, which is physically important.  The time-fractional derivative $\partial_{0+}^\alpha u$ is the Caputo derivative of order $0<\alpha<1$ defined by
\begin{eqnarray*}\label{eq2}
     \partial_{0+}^\alpha u(x,t): =\frac{1}{\Gamma(1-\alpha)}
     \int_0^t \frac{u_\tau(x,\tau)}{(t-\tau)^\alpha }d\tau,
\end{eqnarray*}
which represents the slow diffusion process.

The well-posedness of the time-fractional diffusion system with homogeneous Dirichlet boundary condition corresponding to $\lambda(x)=\infty$ in (\ref{eq1}) has been studied extensively. Sakamoto and Yamamoto \cite{Sakamoto} studied the unique existence of the weak solution and the asymptotic behavior based on the eigenfunction expansions, and the maximal regularity of the solutions in the fractional Sobolev spaces were established in \cite{Gorenflo}. Also, Luchko established a weak maximum principle by a key estimate of the Caputo derivative at an extreme point \cite{Luchko,Luchko1}, from which the uniqueness of a classical solution was also obtained. Based on the weak maximum principle, Liu \cite{Yikan1,Yikan2} improved the maximum principle for fractional diffusion equations with Caputo derivatives. In case of homogeneous Neumann boundary condition corresponding to $\lambda(x)=0$ in (\ref{eq1}), Wei et al proved the existence and uniqueness of a strong solution \cite{TWei1}.
The homogeneous boundary conditions, although they are not fully consistent with the engineering situations, enable the easy mathematical treatments for the governed system such as the direct representations of the solutions in terms of the eigenfunction expansions.

As for inhomogeneous boundary conditions, there are rare works. Recently, Yamamoto \cite{{Yamamoto}} discussed an initial-boundary problem with non-zero Dirichlet boundary value and proved the unique existence of weak solution and the {\it a-priori} estimate for the Caputo fractional derivative in Sobolev spaces. Fujishiro \cite{Fujishiro} considered the approximate controllability with inhomogenous Dirichlet or Neumann boundary condition, and the unique existence and regularity of the solution were shown using the eigenfunction expansions. For the Robin boundary condition with impedance coefficient depending on $(x,t)$, the existence and continuous dependence of the solution in $C(\Omega_T)$ was proved by integral equation method in a bounded domain with Lyapunov boundary \cite{Kemppainen} under the continuous regularity of $(F(x,t), b(x,t),u_0(x))$. In one-dimensional spatial case, Wei and Wang considered the existence and uniqueness of a weak solution with Robin coefficient depending on $t$ \cite{WangJG, TWei3}.

In the practical situations, the requirements of diffusion sources $(F, b, u_0)$ on continuous space are too restrictive, since the measurements of these sources by some average quantities are widely applied in practical situations for smoothing the random errors. On the other hand, the time-depending diffusion process for the contaminants comes from the variable sources with the stable boundary status for most of the practical situations, i.e., the boundary impedance coefficient may be time independent in general. Due to these reasons, we are motivated to consider the solvability of (\ref{eq1}) in some general function space for input sources with relaxed regularities. The well-posedness of direct problem (\ref{eq1}) is important for considering the corresponding inverse problems such as the identification of Robin coefficient $\lambda(x)$, which also motivates us to establish the maximum principles for this system.

In this paper,  we will establish the uniqueness and regularity for either the weak solution or the classical one to (\ref{eq1}) as well as the maximum principles. To this end, instead of using the Hanack inequality, the novel Hopf lemma is established for diffusion operator with time fractional order derivative. More precisely,
In section 2, we present some preliminaries. The well-posedness of the weak and classical solutions for the time-fractional diffusion equation with inhomogenous Robin boundary condition is presented based on the eigenfunction expansions in section 3.
In section 4, we prove the Hopf lemma for time-fractional diffusion operator by constructing an appropriate auxiliary function, generalizing the Hopf lemma for parabolic system with first order time derivative, and then the maximum principles for the initial-boundary problem are established. Finally, some conclusions and future works will be stated in section 5.

\section{Preliminaries}

The representation of the Caputo derivative in terms of the Riemann-Liouville fractional integral for smooth functions is necessary for our analysis.

\begin{defn}
For $y(t)\in L^p(0,T)$ with $1\le p\le +\infty$, the Riemann-Liouville fractional left integral $I_{0+}^\alpha y$ and right integral $I_{T-}^\alpha y$ for $\alpha\in (0,1)$ are defined by
\begin{eqnarray*}
I_{0+}^\alpha y(t):&=&\frac{1}{\Gamma(\alpha)}\int_0^t \frac{y(\tau)}{(t-\tau)^{1-\alpha}}d\tau, \quad 0<t\leq T,\\
I_{T-}^\alpha y(t):&=&\frac{1}{\Gamma(\alpha)}\int_t^T \frac{y(\tau)}{(\tau-t)^{1-\alpha}}d\tau, \quad 0\leq t< T.
\end{eqnarray*}
\end{defn}

It is well-known \cite{Ana} that $I_{0+}^\alpha,I_{T-}^\alpha$ are bounded from $L^p(0,T)$ to $L^p(0,T)$ by
$$\|I_{0+}^\alpha y\|_{L^p(0,T)}, \|I_{T-}^\alpha y\|_{L^p(0,T)}\le \frac{T^\alpha}{\alpha |\Gamma(\alpha)|}\|y\|_{L^p(0,T)}.$$
Moreover, for $p\in (1,\frac{1}{\alpha})$, $I_{0+}^\alpha$ and $I_{T-}^\alpha$ are bounded from $L^p(0,T)$ to $L^q(0,T)$ with $q:=\frac{p}{1-\alpha p}$.

Therefore, for $y\in AC[0,T]$, the space of absolutely continuous functions, its Caputo derivatives are in $L^1(0,T)$ from the relations
\begin{eqnarray*}
\partial _{0+}^{\alpha} y(t):&=&\frac{1}{\Gamma(1-\alpha)}\int_0^t \frac{y'(\tau)}{(t-\tau)^\alpha}d\tau=(I_{0+}^{1-\alpha}y')(t),~~0<t\leq T,\\
\partial _{T-}^{\alpha} y(t):&=&-\frac{1}{\Gamma(1-\alpha)}\int_t^T \frac{y'(\tau)}{(\tau-t)^\alpha}d\tau=-(I_{T-}^{1-\alpha}y')(t),~~0\leq t < T
\end{eqnarray*}
due to $y'\in L^1(0,T)$.

Also, for $0<\alpha\leq 1$, the Riemann-Liouville fractional left derivative and right one are introduced by
\begin{eqnarray*}
D^\alpha_{0+} y(t):=\frac{1}{\Gamma(1-\alpha)}\frac{d}{dt}\int_0^t \frac{y(\tau)}{(t-\tau)^\alpha}d\tau\equiv\frac{d}{dt}(I_{0+}^{1-\alpha}y)(t),~~0<t\leq T,\\
D^\alpha_{T-} y(t):=-\frac{1}{\Gamma(1-\alpha)}\frac{d}{dt}\int_t^T \frac{y(\tau)}{(\tau-t)^\alpha}d\tau\equiv-\frac{d}{dt}(I_{T-}^{1-\alpha}y)(t),~~0\le t< T.\\
\end{eqnarray*}

\begin{prop}\cite{Samko}
Let $\alpha\geq 0$, $p \geq 1$, $q \geq 1$, and $\frac{1}{p}+\frac{1}{q}\leq 1+\alpha$ ($p\neq1$ and $q\neq 1$ in the case $\frac{1}{p}+\frac{1}{q}= 1+\alpha$). If $u(t)\in L^p(0,T)$ and $v(t)\in L^q(0,T)$, then
\begin{eqnarray}\label{De1}
(I_{0+}^\alpha v)\ast u(t)  =v \ast (I_{0+}^\alpha u)(t).
\end{eqnarray}

\end{prop}

%
%

For given $\alpha>0$ and $\beta\in \mathbb{R}$, the Mittag-Leffler function is defined as
\begin{eqnarray*}
 E_{\alpha,\beta}(z):=\sum_{k=0}^{\infty}\frac{z^k}{\alpha k+\beta},~~ z\in \mathbb{C},
\end{eqnarray*}
which is a generalization of the exponential function $e^z$.

%


The following $L^p$ estimate on elliptic system is standard.

\begin{lemma} \label{HtwoM}
Let $w(x)\in H^2(\Omega)$ solve
\begin{eqnarray}
\begin{cases}
-\Delta w=f(x), &x\in\Omega\\
\frac{\partial w}{\partial\nu}+\lambda(x)w(x)=b(x), &x\in\partial\Omega
\end{cases}
\end{eqnarray}
for $0<\lambda_0\le \lambda(x)\in C(\partial\Omega)$. Then we have
$$\norm{w}_{H^2(\Omega)}\le C(\norm{f}_{L^2(\Omega)}+\|\tilde b\|_{H^1(\Omega)})\le C(\norm{f}_{L^2(\Omega)}+\|b\|_{H^{1/2}(\partial\Omega)}),$$
where $\tilde b$ is the $H^1(\Omega)-$extension of $b(x)\in H^{1/2}(\partial\Omega)$.
\end{lemma}

Firstly, we will define the classical solution to direct problem (\ref{eq1}) \cite{{Luchko}}.
\begin{defn}\label{DefReGU11}
The classical solution $u=u(x,t)$ to (\ref{eq1}) is a function  defined in $\Omega_T$ that belongs to the space $C(\overline{\Omega}_T)\cap S_1(\Omega_T)\cap S_2(\Omega_T)$, where
$$W_t^1 (0, T):=\{g \in C^1(0, T ], g'\in L^1(0,T)\}, $$
$$S_1(\Omega_T):=\{F(x,t): F(x,\cdot)\in W_t^1(0,T)\hbox{ for every } x\in\Omega \},$$
$$S_2(\Omega_T):=\{F(x,t): F(\cdot,t)\in C^2(\Omega)\hbox{ for every } t\in (0,T]\}.$$
\end{defn}

Define the Robin eigen-system $\{(\mu_n, \psi_n(x)):n\in \mathbb{N}\}$ of $-\Delta$ in $\Omega$ by
\begin{eqnarray}\label{z7}
 \begin{cases}
 -\Delta \psi_n(x)=\mu_n \psi_n(x), ~~~x\in\Omega\\
 \frac{\partial \psi_n(x)}{\partial\nu} +\lambda(x)\psi_n(x)=0,~~~x\in\partial\Omega
 \end{cases}
\end{eqnarray}
for $0<\lambda(x)\in C(\Omega)$, which satisfies that $0<\mu_1\leq \mu_2 \leq \cdots  \leq$$\mu_n\rightarrow \infty $ as $n\rightarrow \infty$, and the set $\{\psi_n\in H^2(\Omega): n\in \mathbb{N}\}$ forms a normalized orthogonal base of $L^2(\Omega)$ with $\norm{\psi_n}_{L^2(\Omega)}=1$
(The proofs see Appendix ). Then,
we can define the fractional power  operator $(-\Delta)^{\gamma}$  for $\gamma>0$ as
$(-\Delta)^\gamma f:=\sum_{n=1}^{\infty} \mu_n^\gamma (f,\psi_n)\psi_n$
with the domain
$$\mathcal{D}((-\Delta)^\gamma)=\left\{f\in L^2(\Omega): \sum_{n=1}^{\infty} |\mu_n^\gamma (f,\psi_n)|^2<\infty\right\},$$
which is a Hilbert space with the norm
\begin{eqnarray*}
\|f\|_{\mathcal{D}((-\Delta)^\gamma)}= \left(\sum_{n=1}^{\infty} |\mu_n^\gamma (f,\psi_n)|^2\right)^\frac{1}{2}.
\end{eqnarray*}

Since we will construct the explicit solution to our direct problem in terms of the eigen-system formally, we need the higher regularities of eigen-functions from their lower regularities.

\begin{lemma}\label{LemmaPsin}
Assume $\lambda(x)\in H^{k+1/2}(\partial\Omega)$ for given $k\in \{0,1,2,3,\cdots\}$. If the eigen-function $\psi_n\in H^{k+1}(\Omega)$, then $\psi_n\in H^{k+2}(\Omega)$ with the estimate
\begin{eqnarray}\label{ReGuEs1}
 \|\psi_n\|_{H^{k+2}(\Omega)}\leq C\mu_n \|\psi_n\|_{H^{k}(\Omega)}+C\|\psi_n\|_{H^{k+1}(\Omega)},
\end{eqnarray}
where the constant $C$ is independent of $k$.
\end{lemma}

\begin{proof}
We rewrite (\ref{z7}) as
\begin{eqnarray*}\label{z7-1}
 \begin{cases}
 -\Delta \psi_n(x)=\mu_n \psi_n(x):=F_n(x), ~~~x\in\Omega,\\
 \frac{\partial \psi_n}{\partial\nu}=-\lambda(x)\psi_n(x):=G_n(x),~~~x\in\partial\Omega.
 \end{cases}
\end{eqnarray*}
By the regularity estimate for elliptic system with inhomogenous Neumann boundary condition, we have the regularity estimate (Prop.7.5, \cite{Taylor})
$$\|\psi_n\|_{H^{k+2}(\Omega)}\le C\|F_n\|_{H^k(\Omega)}+C\|G_n\|_{H^{k+1/2}(\partial\Omega)}+C\|\psi_n\|_{H^{k+1}(\Omega)},$$
with a constant $C$ independent of $k$. That is,
\begin{eqnarray}\label{CC00-399}
\|\psi_n\|_{H^{k+2}(\Omega)}
\leq C\mu_n \|\psi_n\|_{H^{k}(\Omega)}+C\|\psi_n\|_{H^{k+1/2}(\partial\Omega)}+C\| \psi_n\|_{H^{k+1}(\Omega)}
\end{eqnarray}
from the expressions of $(F_n, G_n)$ and the regularity of $\lambda(x)$. Finally our result follows from the trace theorem (\cite{Lions})
$\|\psi_n\|_{H^{k+1/2}(\partial\Omega)}\leq C \|\psi_n\|_{H^{k+1}(\Omega)}$.

The proof is complete.
\end{proof}

\begin{remark}
For the normalized eigenfunction $\{\psi_n: n=1,2,\cdots\}$ in $L^2$ sense, it is easy to see for  $n=1,2,\cdots$ that
$$\|\psi_n\|_{L^2(\Omega)}=1, \quad \|\psi_n\|_{H^1(\Omega)}\le C(\|\psi_n\|_{L^2(\Omega)}+\|\nabla\psi_n\|_{L^2(\Omega)})\le C\sqrt{\mu_n}.$$
Therefore we can lead to the estimates that
$\|\psi_n\|_{H^k(\Omega)}\le C_k(\sqrt{\mu_n})^k$ for $k=0,1,2,\cdots$
recursively in terms of (\ref{ReGuEs1}).
\end{remark}



\section{The well-posedness of the direct problem }

The main difficulty for showing the regularity of the solution comes from the inhomogeneous Robin boundary condition in (\ref{eq1}). So we firstly consider this case with homogeneous equation and zero initial status, and then apply the known regularity results for inhomogeneous equation.

\begin{lemma}\label{Estimate2}
For $0<\alpha<1$, consider the initial-boundary value problem
\begin{eqnarray}\label{FB4}
\begin{cases}
\partial_{0+}^\alpha z_b(x,t)-\Delta z_b(x,t)=0, &(x,t)\in \Omega_T,\\
\frac{\partial z_b(x,t)}{\partial \nu}+\lambda(x)z_b(x,t)=b(x,t), &(x,t)\in \partial\Omega\times [0,T],\\
z_b(x,0)=0, &x\in\Omega
\end{cases}
\end{eqnarray}
with Robin coefficient $0<\lambda(x)\in C(\partial\Omega)$ and $b\in C([0,T]\times\partial\Omega)$ satisfying the compatibility condition $b(x,0)=0$.

(a) If $b\in C^1([0,T];H^{1/2}(\partial\Omega))$, then there exists a unique solution $z_b\in C([0,T];H^2(\Omega))$
 with the estimate
\begin{eqnarray}\label{FAD}
\|z_b\|_{C([0,T];H^2(\Omega))}\leq  C(\Omega,\alpha,T) \|b\|_{C^1([0,T];H^{1/2}(\partial\Omega))}.
\end{eqnarray}

(b) If $b\in C^2([0,T];H^{1/2}(\partial\Omega))$ satisfying $\partial_t b(x,0)=0$,  then there exists a unique solution $z_b\in C^1([0,T];H^2(\Omega))$
with the estimate
\begin{eqnarray}\label{F12-2}
\|z_b\|_{C^1([0,T];H^2(\Omega))}
 \leq C(\Omega,\alpha,T) \|b\|_{C^2([0,T];H^{1/2}(\partial\Omega))}.
\end{eqnarray}
\end{lemma}

\begin{proof}
We will prove (a). For any fixed $t\in[0,T]$, define $\Lambda_b(\cdot, t)$ from
\begin{eqnarray}\label{elliptic}
\begin{cases}
 -\Delta  \Lambda_b=0, &x\in \Omega,\\
  \frac{\partial \Lambda_b}{\partial \nu}+\lambda(x)\Lambda_b=b(x,t), &x\in\partial\Omega.
\end{cases}
\end{eqnarray}
By the regularity estimate on elliptic problem and Lemma \ref{HtwoM},  there exists a unique solution $\Lambda_b(\cdot,t)\in H^2(\Omega)$ satisfying $\|\Lambda_b(\cdot,t)\|_{H^{2}(\Omega)}\le C_t\|\tilde b(\cdot, t)\|_{H^1(\Omega)}$. Since $\lambda(x)$ is independent of $t$, we also have for $k=0,1,2$ that
\begin{eqnarray}\label{SEP00-1}
\|\partial_t^k \Lambda_b(\cdot,t)\|_{H^2(\Omega)}\leq C \|\partial_t^k \tilde b(\cdot,t)\|_{H^1(\Omega)}\le C\|\partial_t^k b(\cdot,t)\|_{H^{1/2}(\partial\Omega)},\quad t\in[0,T].
\end{eqnarray}
For $b(x,\cdot)\in AC[0,T]$, it follows $\Lambda_b(x,\cdot)\in AC[0,T]$ by (\ref{elliptic}). So we can decompose $z_b=z_c+\Lambda_b$, where $z_c$ satisfies
\begin{eqnarray}\label{FB4-an}
\begin{cases}
\partial_{0+}^\alpha z_c(x,t)-\Delta z_c(x,t)=-\partial_{0+}^\alpha \Lambda_b, &(x,t)\in \Omega_T,\\
\frac{\partial z_c(x,t)}{\partial \nu}+\lambda(x)z_c(x,t)=0, &(x,t)\in \partial\Omega\times [0,T],\\
z_c(x,0)=-\Lambda_b(x,0)=0, &x\in\Omega.
\end{cases}
\end{eqnarray}
Then the solution $z_c(x,t)$  can be represented formally as
\begin{eqnarray*}
z_c(x,t)=-\sum_{n=1}^{\infty} \int_0^t \left(\partial_{0+}^\alpha \Lambda_b(\cdot,t-\tau),\psi_n(\cdot)\right)\tau^{\alpha-1}E_{\alpha,\alpha}(-\mu_n \tau^\alpha) d\tau  \;\psi_n(x).
\end{eqnarray*}
Using the regularity of $\Lambda_b(x,t)$ and (\ref{De1}), the above representation says
\begin{eqnarray}\label{ZC-1}
& &z_c(x,t)\nonumber\\&=&-\sum_{n=1}^{\infty} \int_0^t \left( I_{0+}^{1-\alpha}\frac{d}{dt} \Lambda_b(\cdot,t-\tau),\psi_n(\cdot)\right)\tau^{\alpha-1}E_{\alpha,\alpha}(-\mu_n \tau^\alpha) d\tau  \;\psi_n(x) \nonumber\\
&=& -\sum_{n=1}^{\infty} \int_0^t \left( \frac{d}{d t} \Lambda_b(\cdot,t-\tau),\psi_n(\cdot)\right)\cdot I_{0+}^{1-\alpha}\left(\tau^{\alpha-1}E_{\alpha,\alpha}(-\mu_n \tau^\alpha) \right)d\tau  \;\psi_n(x) \nonumber\\
&=& -\sum_{n=1}^{\infty} \int_0^t \left( \partial_t \Lambda_b(\cdot,t-\tau),\psi_n(\cdot)\right)\cdot E_{\alpha,1}(-\mu_n \tau^\alpha) d\tau  \;\psi_n(x),
\end{eqnarray}
which leads to
\begin{eqnarray*}\label{July1}
\|z_c(\cdot,t)\|^2_{L^2(\Omega)}&\le&
\sum_{n=1}^{\infty}\int_0^t \left(\partial_t\Lambda_b(\cdot,t-\tau),\psi_n(\cdot)\right)^2d\tau\;
\int_0^t \left|E_{\alpha,1}(-\mu_n \tau^\alpha)\right|^2d\tau
\nonumber\\&\le&
\sum_{n=1}^{\infty}\int_0^t \left(\partial_\tau\Lambda_b(\cdot,\tau),\psi_n(\cdot)\right)^2d\tau\;
\int_0^t\left(\frac{C}{1+\mu_n \tau^\alpha}\right)^2d\tau
\nonumber\\&\le&
C^2t\int_0^t \sum_{n=1}^{\infty}\left(\partial_\tau\Lambda_b(\cdot,\tau),\psi_n(\cdot)\right)^2d\tau
\nonumber\\&=&
C^2t\int_0^t \|\partial_\tau\Lambda_b(\cdot,\tau)\|^2_{L^2(\Omega)}d\tau\le C^2t\|\partial_t\Lambda_b\|^2_{L^2((0,T)\times \Omega)},
\end{eqnarray*}
that is, $z_c\in C([0,T], L^2(\Omega))$ satisfies
\begin{equation}
\|z_c\|_{C([0,T], L^2(\Omega))}\le C\|\partial_tb\|_{L^2((0,T); H^{1/2}(\partial\Omega))}.
\end{equation}
Therefore we have from the decomposition that
\begin{eqnarray}\label{Czb0-0}
\nonumber \|z_b\|_{C([0,T];L^2(\Omega))}
&\leq& \|\Lambda_b\|_{C([0,T];L^2(\Omega))}+\|z_c\|_{C([0,T];L^2(\Omega))}
\\  \nonumber
&\leq&
C(\|b\|_{C([0,T];H^{1/2}(\partial\Omega))}+\|\partial_tb\|_{C([0,T];H^{1/2}(\partial\Omega))})
\\
&\leq&
C\|b\|_{C^1([0,T];H^{1/2}(\partial\Omega))}.
\end{eqnarray}
Now we estimate $\|\Delta z_c(\cdot,t)\|_{L^2(\Omega)}$. Taking Laplacian operation to (\ref{ZC-1}) says
\begin{eqnarray}\label{July2}
& & \|\Delta z_c(\cdot,t)\|^2_{L^2(\Omega)}
\nonumber\\&=&
\sum_{n=1}^{\infty} \mu_n^2\left(\int_0^t \left(  \partial_t\Lambda_b(\cdot,t-\tau),\psi_n(\cdot)\right)\cdot E_{\alpha,1}(-\mu_n \tau^\alpha) d\tau \right)^2 \nonumber\\
&\le&\sum_{n=1}^{\infty}\mu_n^{2} \int_0^t \left|\left(
\partial_t\Lambda_b(\cdot,t-\tau),\psi_n(\cdot)\right)\right|^2d\tau \int_0^t\left|E_{\alpha,1}(-\mu_n \tau^\alpha)\right|^2 d\tau \nonumber\\
&=& \sum_{n=1}^{\infty}\mu_n^{2\theta} \int_0^t\left|\partial_t \Lambda_b(\cdot,t-\tau),\psi_n(\cdot)\right|^2d\tau \cdot\int_0^t \left|\mu_n^{1-\theta}E_{\alpha,1}(-\mu_n \tau^\alpha)\right|^2  d\tau\quad\quad
\end{eqnarray}
for $1/2\le \theta<1$. By the property of $E_{\alpha,1}(\cdot)$, we have
\begin{eqnarray*}
\left|\mu_n^{1-\theta}E_{\alpha,1}(-\mu_n \tau^\alpha)\right|\leq \mu_n^{1-\theta} \frac{C}{1+\mu_n \tau^\alpha}=C \frac{(\mu_n \tau^\alpha)^{1-\theta}}{1+\mu_n \tau^\alpha}
\cdot \tau^{\alpha(\theta-1)}\leq C\tau^{\alpha(\theta-1)}.
\end{eqnarray*}
Therefore,
\begin{eqnarray}\label{July3}
& & \|\Delta z_c(\cdot,t)\|^2_{L^2(\Omega)}
\nonumber\\&\leq&
\sum_{n=1}^{\infty}\mu_n^{2\theta}\int_0^t \left(\partial_t \Lambda_b(\cdot,t-\tau),\psi_n(\cdot)\right)^2 d\tau \frac{1}{1+2\alpha(\theta-1)}t^{1+2\alpha(\theta-1)}\nonumber\\
 &=& Ct^{1+2\alpha(\theta-1)}\int_0^t
 \left\|\partial_t\Lambda_b(\cdot,t-\tau)\right\|^2_{D(-(\Delta)^\theta)} d\tau \nonumber\\
 &\le& Ct^{1+2\alpha(\theta-1)}\|\partial_t\Lambda_b(\cdot,\cdot)\|^2_{L^2(0,T;D(-(\Delta)^\theta)}.
 \end{eqnarray}
Especially, it follows from $\|\cdot\|_{D(-\Delta)^{\theta}}\sim \|\cdot\|_{H^{2\theta}}$
for $\theta=\frac{1}{2}$ that
\begin{eqnarray}\label{July3-1}
\nonumber \|\Delta z_c(\cdot,t)\|^2_{L^2(\Omega)}\le Ct^{1-\alpha}\|\partial_t\Lambda_b(\cdot,\cdot)\|^2_{L^2(0,T;H^1(\Omega))}\le Ct^{1-\alpha}
\|\partial_t b\|^2_{L^2(0,T; H^{1/2}(\partial\Omega))},
\end{eqnarray}
which leads to $\lim_{t\rightarrow 0}\| \Delta z_c(\cdot,t)\|_{L^2(\Omega)}=0$, i.e.,  $\Delta z_c\in C([0,T], L^2(\Omega))$ satisfies
\begin{equation}
\|\Delta z_c\|_{C([0,T], L^2(\Omega))}\le C\|\partial_tb\|_{L^2(0,T; H^{1/2}(\partial\Omega))}.
\end{equation}
Finally we have from the decomposition that
 \begin{eqnarray}\label{Mm10}
  \nonumber \|\Delta z_b(\cdot,t)\|_{L^2(\Omega)}^2
 &\leq &\|\Delta z_c(\cdot,t)\|_{L^2(\Omega)}^2+\|\Delta \Lambda_b (\cdot,t)\|_{L^2(\Omega)}^2\\ \nonumber
  &\leq &C\left(t^{1-\alpha}
\|\partial_t b\|^2_{L^2(0,T; H^{1/2}(\partial\Omega))}+\| \Lambda_b (\cdot,t)\|_{H^2(\Omega)}^2\right)\\  \nonumber
 &\leq &C\left(t^{1-\alpha}
\|\partial_t b\|^2_{L^2(0,T; H^{1/2}(\partial\Omega))}+\| b (\cdot,t)\|_{H^{1/2}(\partial\Omega)}^2\right).
 \end{eqnarray}
For any fixed $t>0$, we  utilize Lemma \ref{HtwoM} to get
\begin{eqnarray}\label{July4}
\|z_b(\cdot,t)\|_{H^2(\Omega)}
&\leq& \|\Delta z_b(\cdot,t)\|_{L^2(\Omega)}+\|b(\cdot,t)\|_{H^{1/2}(\partial\Omega)}\nonumber \\
&\leq& C\left( t^{(1-\alpha)/2} \left\| \partial_t b  \right\|_{L^2(0,T;H^{1/2}(\partial\Omega))}+\|b(\cdot,t)\|_{H^{1/2}(\partial\Omega)}\right).\;\qquad
\end{eqnarray}
Thus $z_b\in C([0,T];H^2(\Omega))$ for $b\in C^1([0,T], H^{1/2}(\partial\Omega))$ with the estimate
\begin{eqnarray}\label{July4-11}
\|z_b\|_{C([0,T];H^2(\Omega))}
&\leq& C\left(  \left\| \partial_t b  \right\|_{C([0,T];H^{1/2}(\partial\Omega))}+\|b\|_{C([0,T];H^{1/2}(\partial\Omega))}\right)
\nonumber \\&\leq&
C \left\|b\right\|_{C^1([0,T];H^{1/2}(\partial\Omega))},
\end{eqnarray}
where $C=C(\Omega,T,\alpha)$, which verifies (\ref{FAD}).

The result (b) is essentially the same as (a) for $b$ with higher regularity with respect to $t$.

The proof is complete.
\end{proof}

Now we can show the regularity of the solution to (\ref{eq1}).

\begin{theorem}\label{Weak}
For $0<\alpha<1$ and $0<\lambda_0\le \lambda(x)\in C(\partial\Omega)$, we assume the
compatibility condition
$$\frac{\partial u_0(x)}{\partial \nu(x)}+\lambda(x)u_0(x)=b(x,0),\quad x\in \partial\Omega$$
and $(f,g)\in L^2(\Omega)\times L^\infty (0,T)$ and $b\in C^1([0,T];H^{1/2}(\partial\Omega))$.
Then

(i) For $u_0\in L^2(\Omega)$, there exists a unique weak solution $u\in C([0,T];L^2(\Omega))$.
Moreover, there exists a constant $\widetilde{C}_1(\Omega,\alpha,T)>0$   such that
\begin{eqnarray}\label{AA1}
\nonumber \|u\|_{C([0,T];L^2(\Omega))}
 &\leq& \widetilde{C}_1 (\|u_0\|_{L^2(\Omega)}+\|g\|_{L^\infty(0,T)}\|f\|_{L^2(\Omega)}+\\
 &&\hskip 0.6cm\|b\|_{C^1([0,T];H^{1/2}(\partial\Omega))}).
\end{eqnarray}

(ii) For $u_0(x), f(x)\in H^2(\Omega)$, there exists a unique weak solution $u\in C([0,T];H^2(\Omega))$.  Moreover, there exists a constant $\widetilde{C}_2(\Omega,\alpha,T)>0$ such that
\begin{eqnarray}\label{AA2}
\nonumber \|u\|_{C([0,T];H^2(\Omega))}
 &\leq& \widetilde{C}_2 (\|u_0\|_{H^2(\Omega)}+\|g\|_{L^\infty(0,T)}\|f\|_{H^2(\Omega)}+\\
 &&\hskip 0.6cm \|b\|_{C^1([0,T];H^{1/2}(\partial\Omega))}).
\end{eqnarray}
\end{theorem}

\begin{proof}
(i) By the linear superposition principle and Lemma \ref{Estimate2}, the solution of (\ref{eq1}) has the formal representation
\begin{eqnarray}\label{AA3}
u(x,t)&=& \sum_{n=1}^{\infty} u_{0n} \; E_{\alpha,1}(-\mu_n t^\alpha) \psi_n(x)+\nonumber\\
& &\sum_{n=1}^{\infty} f_n\int_0^t g(t-\tau)\tau^{\alpha-1}E_{\alpha,\alpha}(-\mu_n \tau^\alpha)  d\tau \psi_n(x)+u_3(x,t)\nonumber\\
&=:& u_1(x,t)+u_2(x,t)+u_3(x,t),
\end{eqnarray}
where $u_3(x,t)$ solves (\ref{FB4}). We estimate each term respectively. Firstly, we have
\begin{eqnarray}\label{AA4}
\|u_1(\cdot,t)\|_{L^2(\Omega)}^2= \sum_{n=1}^{\infty} u^2_{0n} \; |E_{\alpha,1}(-\mu_n t^\alpha)|^2
\leq C \sum_{n=1}^{\infty} u^2_{0n}=C \|u_0\|^2_{L^2(\Omega)}.
\end{eqnarray}
Introducing
\begin{eqnarray*}\label{AA5}
g_n(t):=\int_0^t g(t-\tau)\tau^{\alpha-1}E_{\alpha,\alpha}(-\mu_n \tau^\alpha)  d\tau, \quad n\in \mathbb{N}.
\end{eqnarray*}
Then it follows for $t\in [0,T]$ that
\begin{eqnarray*}\label{AA5-0}
|g_n(t)|\leq \|g\|_{L^\infty(0,T)} t^\alpha /\Gamma(1+\alpha)\leq \|g\|_{L^\infty(0,T)} T^\alpha /\Gamma(1+\alpha),
\end{eqnarray*}
\begin{eqnarray*}\label{AA5-1}
|g_n(t)|\leq \|g\|_{L^\infty(0,T)}(1-E_{\alpha,1}(-\mu_n t^\alpha))/\mu_n \leq  \|g\|_{L^\infty(0,T)}/\mu_n,
\end{eqnarray*}
which yields
\begin{eqnarray}\label{AA6}
\|u_2(\cdot,t)\|_{L^2(\Omega)}^2= \sum_{n=1}^{\infty} f^2_{n} \; |g_n(t)|^2
\leq \|g\|^2_{L^\infty(0,T)}\frac{t^{2\alpha}}{\Gamma^2(1+\alpha)}
\|f\|^2_{L^2(\Omega)}.
\end{eqnarray}
We can see that
$\lim_{t\rightarrow 0} \|u_2(\cdot,t)\|_{L^2(\Omega)}=0$ due to $u_2(x,0)=0$ and therefore
$\|u_2(\cdot,t)\|_{L^2(\Omega)}\in C[0,T]$.
As for $u_3(x,t)$, we have from (\ref{Czb0-0}) that
\begin{eqnarray}\label{AA7}
\|u_3(\cdot,t)\|_{L^2(\Omega)}^2 \leq  C \|b\|_{C^1([0,T];H^{1/2}(\partial\Omega))}^2
\end{eqnarray}
and $\lim_{t\rightarrow 0} \|u_3(\cdot,t)\|_{L^2(\Omega)}=0$ by the estimates on $\Lambda_b(\cdot,t)$ and $z_c(\cdot,t)$.
From (\ref{AA4}), (\ref{AA6}) and (\ref{AA7}), we obtain
\begin{eqnarray*}
\|u(\cdot,t)\|_{L^2(\Omega)} \leq  \widetilde{C}_1(\|u_0\|_{L^2(\Omega)}+\|g\|_{L^\infty(0,T)}\|f\|_{L^2(\Omega)}+ \|b\|_{C^1([0,T];H^{1/2}(\partial\Omega))}),
\end{eqnarray*}
where $\widetilde{C}_1$ depends on $(\Omega,\alpha,T)$.
On the other hand, we have
\begin{eqnarray}\label{liu3}
& &\|u(\cdot,t)-u_0(\cdot)\|_{L^2(\Omega)}^2
\nonumber\\&\leq& \sum_{n=1}^{\infty} u_{0n}^2 \; (E_{\alpha,1}(-\mu_n t^\alpha)-1)^2+\|u_2(\cdot,t)\|_{L^2(\Omega)}^2+\|u_3(\cdot,t)\|_{L^2(\Omega)}^2.\quad
\end{eqnarray}
Noticing
\begin{eqnarray*}
& &\sum_{n=1}^{\infty} u_{0n}^2 \; (E_{\alpha,1}(-\mu_n t^\alpha)-1)^2\nonumber\\&=&
\sum_{n=1}^{M} u_{0n}^2 \; (E_{\alpha,1}(-\mu_n t^\alpha)-1)^2+
\sum_{n=M+1}^{\infty} u_{0n}^2 \; (E_{\alpha,1}(-\mu_n t^\alpha)-1)^2
\nonumber\\&\le&
\sum_{n=1}^{M} u_{0n}^2 \; (E_{\alpha,1}(-\mu_n t^\alpha)-1)^2+
\sum_{n=M+1}^{\infty} u_{0n}^2
\end{eqnarray*}
due to $0<E_{\alpha,1}(-\mu_n t^\alpha)\le 1$
for $t\in [0,T]$, the first term  in (\ref{liu3}) tends to zero as $t \to 0$, and therefore we have the continuity of $u_1$. Finally we have
\begin{eqnarray*}
\lim_{t\rightarrow 0} \|u(\cdot,t)-u_0\|_{L^2(\Omega)}=0.
\end{eqnarray*}
So $u\in C([0,T];L^2(\Omega))$.

(ii) We will verify $u \in C([0,T];H^2(\Omega))$ and (\ref{AA2}).
Firstly, we have
\begin{eqnarray}\label{AA9}
\|\Delta u_1(\cdot,t)\|_{L^2(\Omega)}^2
=\sum_{n=1}^{\infty} u^2_{0n}\mu_n ^2 \; |E_{\alpha,1}(-\mu_n t^\alpha)|^2
\leq C \|u_0\|^2_{H^2(\Omega)}
\end{eqnarray}
by $|E_{\alpha,1}(-\mu_n t^\alpha)|\le C$ and
\begin{eqnarray}\label{AA10}
\|\Delta u_2(\cdot,t)\|_{L^2(\Omega)}^2
=\sum_{n=1}^{\infty} \mu_n ^2 f^2_{n} \; |g_n(t)|^2
\leq\|g\|^2_{L^\infty(0,T)}\frac{t^{2\alpha}}{\Gamma^2(1+\alpha)}
 \|f\|^2_{H^2(\Omega)}
\end{eqnarray}
using (\ref{AA5-0}), we can also see that $\lim_{t\rightarrow 0}\| \Delta u_2(\cdot,t)\|_{L^2(\Omega)}=0$.
As for $\Delta u_3(x,t)$, we have
\begin{eqnarray}\label{AA-LIU}
 \|\Delta u_3(\cdot,t)\|_{L^2(\Omega)}^2
\leq  C\left(t^{1-\alpha}
\|\partial_t b\|^2_{L^2(0,T; H^{1/2}(\partial\Omega))}+\| b (\cdot,t)\|_{H^{1/2}(\partial\Omega)}^2\right)
\end{eqnarray}
from (\ref{Mm10}).
Combining (\ref{AA9})-(\ref{AA-LIU}), we have
\begin{eqnarray*}
  \|\Delta u(\cdot,t)\|_{L^2(\Omega)}
&\leq& C(\|u_0\|_{H^2(\Omega)}+\|g\|_{L^\infty(0,T)}\frac{t^{\alpha}}{\Gamma(1+\alpha)}
 \|f\|_{H^2(\Omega)}+\\
&&\hskip 0.5cm t^{(1-\alpha)/2}
\|\partial_t b\|_{L^2(0,T; H^{1/2}(\partial\Omega))}+\| b (\cdot,t)\|_{H^{1/2}(\partial\Omega)}).
\end{eqnarray*}
Thus we have from Lemma \ref{HtwoM} that
\begin{eqnarray*}
\|u(\cdot,t)\|_{H^2(\Omega)}
&\leq&  C\left(  \|\Delta u(\cdot,t)\|_{L^2(\Omega)}+\|b(\cdot,t)\|_{H^{1/2}(\partial\Omega)}\right)\\
&\leq& \widetilde{C}_2(\|u_0\|_{H^2(\Omega)}+\|g\|_{L^\infty(0,T)}\frac{t^{\alpha}}{\Gamma(1+\alpha)}
 \|f\|_{H^2(\Omega)}+\\
&&\hskip 0.6cm t^{(1-\alpha)/2}
\|\partial_t b\|_{L^2(0,T; H^{1/2}(\partial\Omega))}+\|b\|_{C([0,T],H^{1/2}(\partial\Omega))}).
\end{eqnarray*}
Moreover, we have for $t\ge 0$ that
\begin{eqnarray}\label{eq329}
& &\|u(\cdot,t)-u_0\|_{H^2(\Omega)}
\nonumber\\&\leq& \|\Delta u(\cdot,t)-\Delta u_0\|_{L^2(\Omega)}+ \|b(\cdot,t)-(\frac{\partial u_0}{\partial \nu}+\lambda(x)u_0(x))\|_{H^{1/2}(\partial\Omega)}\quad
\end{eqnarray}
using Lemma \ref{HtwoM} again.
However, it follows from
\begin{eqnarray*}
\Delta u(\cdot,t)
= -\sum_{n=1}^{\infty} u_{0n}\mu_n  \; E_{\alpha,1}(-\mu_n t^\alpha) \psi_n(x)+\Delta u_2(\cdot,t)
+\Delta z_c(\cdot,t)+\Delta \Lambda_b(\cdot,t),\\
\lim_{t\rightarrow 0} \|\Delta z_c(\cdot,t)\|_{L^2(\Omega)}=\lim_{t\rightarrow 0} \|\Delta u_2(\cdot,t)\|_{L^2(\Omega)}=0,\\
\Delta u_0(x)=-\sum_{n=1}^{\infty} u_{0n}\mu_n  \psi_n(x),\quad \Delta \Lambda_b(x,t)=0
\end{eqnarray*}
that
\begin{eqnarray}\label{liu4}
\lim_{t\rightarrow 0} \left\|\Delta u(\cdot,t)-\Delta u_0\right\|_{L^2(\Omega)}^2
\le C\lim_{t\rightarrow 0} \sum_{n=1}^{\infty} u_{0n}^2\mu_n^2  \; (E_{\alpha,1}(-\mu_n t^\alpha)-1)^2=0
\end{eqnarray}
due to $u_0\in H^2(\Omega)$ by the similar arguments.
Considering the compatibility condition,  it follows from (\ref{eq329}) and (\ref{liu4}) that
$$\lim_{t\to 0}\|u(\cdot,t)-u_0(\cdot)\|_{H^2(\Omega)}=0.$$
Therefore, we  obtain that $u \in C([0,T];H^2(\Omega))$.

Finally, we prove the uniqueness of the solution to to problem (\ref{eq1}), i.e., for $u_0=fg=b\equiv 0$, (\ref{eq1}) has only zero solution $u$. Since $\{(\psi_n(x),\mu_n): n\in \mathbb{N}\}$ is the eigen-system,
taking the inner product of (\ref{eq1}) with respect to $\psi_n(x)$ and setting $u_n(t):=(u(\cdot,t),\psi_n)$, we obtain
\begin{eqnarray*}
\begin{cases}
\partial_{0+}^\alpha u_n(t)=-\mu_n u_n(t),~~t\in(0,T],\\
u_n(0)=0.
\end{cases}
\end{eqnarray*}
Due to the well-posedness of the ordinary fractional differential equation
(Ch.3, \cite{Podlubny}), we obtain that $u_n(t)=0,\; n\in \mathbb{N}$.  Since $\{\psi_n(x)\}_{n\in \mathbb{N}}$ is an orthonormal basis in $L^2(\Omega)$, we have $u=0$ in $\Omega\times (0,T]$.

The proof is complete.
\end{proof}

The weak solution, as analyzed above, is always taking some average process for the spatial variable in $L^2$ or $H^2$, which is not applicable for the pointwise measurement. If we enhance the regularity of the input data, we can also have the pointwise estimate on the function $u$ as well as its fractional derivative.

\begin{theorem}\label{Classical}
Assume that $\lambda(x)>0$ is smooth enough, let
$$u_0\in \mathcal{D}((-\Delta)^{3}), \; g\in C^1[0,T],\;f \in \mathcal{D}((-\Delta)^{2}),$$ and
$b\in C^2([0,T]; H^{9/2}(\partial\Omega)) $
Then, there exists a unique classical solution $u\in C(\overline{\Omega}_T)\bigcap S_1(\Omega_T)\bigcap S_2(\Omega_T)$ for (\ref{eq1}). Moreover, there exists a constant $C>0$ depending on $(\Omega,\alpha,T)$ such that the solution has the estimate
\begin{eqnarray}\label{AA00-0}
 & &\| u\|_{C([0,T]\times\overline{\Omega})} \nonumber\\
&\leq& C\left(\|u_0\|_{\mathcal{D}((-\Delta)^2)}+\|g\|_{C[0,T]} \|f\|_{ \mathcal{D}((-\Delta)^{})}+\|b\|_{C^1([0,T]; H^{5/2}(\partial\Omega))}\right).\quad
\end{eqnarray}

\end{theorem}

\begin{proof}
Since (\ref{eq1}) has a formal solution $u$ given by (\ref{AA3}),  it is enough to show $u_i\in C(\overline{\Omega}_T)$ for $i=1,2,3$ for verifying $u\in C(\overline{\Omega}_T)$. To this end, we show that the series defining $u_i$ are uniformly convergent in $\overline \Omega_T$.
By Lemma \ref{LemmaPsin}, we have
\begin{eqnarray}\label{AA00-3}
 \|\psi_n\|_{C(\overline{\Omega})}\leq \|\psi_n\|_{H^{2k}(\Omega)}
\leq C \mu_n^k\le Cn^{2k/d}, ~k=0,1,2,\cdots
\end{eqnarray}
for any $k>d/4$ due to $\mu_n = O(n^{2/d})$. We take $k,m>d/4$ such that $\gamma:=m+k>d/2$.

For the series defining $u_1$ and $u_2$, we have the estimate
\begin{eqnarray}\label{U1-LIU}
& &\sum_{n=1}^{\infty} |u_{0n} \; E_{\alpha,1}(-\mu_n t^\alpha)| \; |\psi_n(x)|\nonumber\\
&\leq& \sum_{n=1}^{\infty} |u_{0n}| \;\|\psi_n\|_{C(\overline{\Omega})}\nonumber\\&\le& C \left( \sum_{n=1}^{\infty}\frac{1}{\mu_n^{2m}}\right)^{1/2} \; \left( \sum_{n=1}^{\infty}\mu_n^{2(k+m)}u_{0n} ^2\right)^{1/2}\le C\|u_0\|_{\mathcal{D}((-\Delta)^\gamma)}
\end{eqnarray}
and
\begin{eqnarray}\label{U2-LIU}
& &\sum_{n=1}^{\infty} |f_n \; g_n(t)| \; |\psi_n(x)|
\nonumber\\&\leq& \|g\|_{C[0,T]} \sum_{n=1}^{\infty}\frac{1}{\mu_n}|f_n|\; \|\psi_n\|_{C(\overline{\Omega})}
\nonumber\\&\leq&
C\|g\|_{C[0,T]} \left( \sum_{n=1}^{\infty}\frac{1}{\mu_n^{2m}}\right)^{1/2}\;
\left( \sum_{n=1}^{\infty}\mu_n^{2(k+m-1)}f_{n} ^2\right)^{1/2}
\nonumber\\&\leq&
C\|g\|_{C[0,T]}
\|f\|_{\mathcal{D}((-\Delta)^{\gamma-1})}.
\end{eqnarray}
As for $u_3(x,t):=z_b(x,t)\equiv z_c(x,t)+\Lambda_b(\cdot,t)$ defining in (\ref{AA3}) solving (\ref{FB4}), since $\Lambda_b(\cdot,\cdot)\in C(\overline\Omega_T)$ with
$$\|\Lambda_b(\cdot, t)\|_{C(\overline\Omega)}\le C \|\Lambda_b(\cdot,t)\|_{H^2(\Omega)}\le C\|b(\cdot,t)\|_{H^{1/2}(\partial\Omega)}$$
for $t\in [0,T]$ by Lemma \ref{HtwoM},  let us show $z_c\in C(\overline\Omega_T)$. For the series defining $z_c(x,t)$, we also have
\begin{eqnarray}\label{U3-LIU}
& &\sum_{n=1}^{\infty} \;  \int_0^t \left| (\partial_t\Lambda_b(\cdot,t-\tau),\psi_n(\cdot))E_{\alpha,1}(-\mu_n \tau^\alpha) \right| d\tau \; |\psi_n(x)|\nonumber\\
&\leq&\sum_{n=1}^{\infty} \|\psi_n\|_{C(\overline{\Omega})}\left(\int_0^t (\partial_t\Lambda_b(\cdot,t-\tau),\psi_n(\cdot))^2 d\tau\right)^{1/2}
\nonumber\\
&\le&C\left( \sum_{n=1}^{\infty}\frac{1}{\mu_n^{2m}}\right)^{1/2} \;
\left( \sum_{n=1}^{\infty}\mu_n^{2(k+m)}\int_0^t(\partial_\tau\Lambda_b(\cdot,\tau),\psi_n(\cdot))^2d\tau\right)^{1/2}
\nonumber\\&\le&
C\left( \sum_{n=1}^{\infty}\frac{1}{\mu_n^{2m}}\right)^{1/2} \;\left( \sum_{n=1}^{\infty}\mu_n^{2(k+m)}\int_0^T(\partial_\tau\Lambda_b(\cdot,\tau),\psi_n(\cdot))^2d\tau\right)^{1/2}
\nonumber\\
&\le&C\|\partial_t\Lambda_{ b}\|_{L^2(0,T;D((-\Delta)^\gamma))}.
\end{eqnarray}
Therefore, for $(u_0,f,g,b)$ satisfying the specified regularities, the three series defining $u_1,u_2,z_c$ are absolutely uniformly convergent in $\overline \Omega_T$ in terms of (\ref{U1-LIU})-(\ref{U3-LIU}), which ensure that $u_1,u_2, u_3$, and consequently $u(x,t)=u_1+u_2+u_3$, are continuous functions in  $\overline \Omega_T$.

Since we consider the cases $d=1,2,3$, we take $\gamma=2$.
So we have
\begin{eqnarray}\label{SEP20-1}
\|\partial_t\Lambda_b\|_{L^2(0,T; \mathcal{D}((-\Delta)^{\gamma})}^2&=&
\int_0^T \sum_{n} \mu_n^{4}(\partial_t\Lambda_b(\cdot,t),\psi_n)^2dt
\nonumber\\
&\leq& C \|\partial_t\Lambda_b\|^2_{L^2(0,T;H^4(\Omega))}.
\end{eqnarray}
By the estimate (Theorem 15.2, \cite{Agmon}) on (\ref{elliptic}), we have
for any integer $l\ge 2$ that
$$\|\partial_t^k\Lambda_b(\cdot,t)\|_{H^l(\Omega)}\le C_l(\|\partial_t^kb(\cdot,t)-\lambda(\cdot)\partial_t^k\Lambda_b(\cdot,t)\|_{H^{l-3/2}(\partial\Omega)}+
\|\partial_t^k\Lambda_b(\cdot,t)\|_{L^2(\Omega)}).$$
Since $\lambda\in H^{l-3/2}(\partial\Omega)$, $\|\partial_t^k\Lambda_b(\cdot,t)\|_{H^{l-3/2}(\partial\Omega)}\le C\|\partial_t^k\Lambda_b(\cdot,t)\|_{H^{l-1}(\Omega)}$, we have
\begin{eqnarray*}
& &\|\partial_t^k\Lambda_b(\cdot,t)\|_{H^l(\Omega)}
\nonumber\\&\le& C_l(\|\partial_t^kb(\cdot,t)\|_{H^{l-3/2}(\partial\Omega)}+
\|\partial_t^k\Lambda_b(\cdot,t)\|_{H^{l-1}(\Omega)}+
\|\partial_t^k\Lambda_b(\cdot,t)\|_{L^2(\Omega)}),
\end{eqnarray*}
which leads to for $l=2,3,\cdots$ that
\begin{equation}
\|\partial_t^k\Lambda_b(\cdot,t)\|_{H^l(\Omega)}\le C_l(\|\partial_t^kb(\cdot,t)\|_{H^{l-3/2}(\partial\Omega)}+\|\partial_t^k\Lambda_b(\cdot,t)\|_{L^2(\Omega)}).
\end{equation}

On the other hand, for $\lambda(x)>0$ on $\partial\Omega$, by the homogenization of the boundary condition, we have the estimate that
$$\|\partial_t^k\Lambda_b(\cdot,t)\|_{L^2(\Omega)}\le \|\partial_t^k\Lambda_b(\cdot,t)\|_{H^2(\Omega)}\le
C\|\partial_t^kb(\cdot,t)\|_{H^{1/2}(\partial\Omega)}$$
from (\ref{SEP00-1}).
So we have for $l=2,3,\cdots$ that
\begin{equation}\label{liu-3333}
\|\partial_t^k\Lambda_b(\cdot,t)\|_{H^l(\Omega)}\le C_l\|\partial_t^kb(\cdot,t)\|_{H^{l-3/2}(\partial\Omega)},
\end{equation}
which leads to
\begin{eqnarray}\label{SEP20-3}
\|\partial_t\Lambda_b\|_{L^2(0,T; \mathcal{D}((-\Delta)^{\gamma})}^2
\leq C \|\partial_t\Lambda_b\|^2_{L^2(0,T;H^4(\Omega))}
\leq C \|\partial_t b\|^2_{L^2(0,T;H^{5/2}(\partial\Omega))}.
\end{eqnarray}
Therefore
\begin{eqnarray}\label{SEP20-2}
\nonumber \| u\|_{C([0,T]\times \overline{\Omega})}
&\leq& C\|u_0\|_{\mathcal{D}((-\Delta)^2)}+C\|g\|_{C[0,T]}
\|f\|_{\mathcal{D}((-\Delta)^{})}+\\
&&C\|\partial_t b\|_{L^2(0,T;H^{5/2}(\partial\Omega))}+
C\|b\|_{C([0,T];H^{1/2}(\partial\Omega))}.
\end{eqnarray}

Next we verify that $u(x,\cdot)\in W_t^1(0,T)$ and $u(\cdot,t)\in C^2(\Omega)$.
For fixed $x\in \Omega$, we firstly verify $u(x,\cdot)\in W_t^1(0,T)$. Similarly, it is enough to verify $u_i(x,\cdot)\in W_t^1(0,T)$ for $i=1,2,3$, i.e., we need to show that the series defines $u_i(x,\cdot)$ with the regularities $u_i(x,\cdot)\in C^1(0,T]$ and $\partial_t u_i(x,\cdot)\in L^1(0,T)$.

We will prove $u_i(x,\cdot)\in C^1(0,T]$. For the series defining $\partial_t u_1$ and $\partial_t u_2$, we have the estimates
\begin{eqnarray}\label{U1t}
\nonumber && \sum_{n=1}^{\infty} \mu_n |u_{0n} \; t^{\alpha-1} E_{\alpha,\alpha}(-\mu_n t^\alpha) | \; |\psi_n(x)|\\ \nonumber
&\leq& \frac{1}{\Gamma(\alpha)} \sum_{n=1}^{\infty} \mu_n \; |u_{0n}| \; t^{\alpha-1}\;\|\psi_n\|_{C(\overline{\Omega})}\\ \nonumber
&\leq& C\frac{1}{\Gamma(\alpha)} \left( \sum_{n=1}^{\infty}\frac{1}{\mu_n^{2m}}\right)^{1/2} \; \left( \sum_{n=1}^{\infty}\mu_n^{2(k+m+1)}u_{0n} ^2\right)^{1/2} t^{\alpha-1}\\
&\leq& C \|u_0\|_{\mathcal{D}((-\Delta)^{\gamma+1})}t^{\alpha-1}
  \end{eqnarray}
  and
 \begin{eqnarray}\label{U2t}
\nonumber && \sum_{n=1}^{\infty}  |g(0) f_n\;  t^{\alpha-1}E_{\alpha,\alpha}(-\mu_n t^\alpha)|  \; |\psi_n(x)|+\\ \nonumber
&&\sum_{n=1}^{\infty}  |f_n|\; \left|\int_0^t \partial_t g(t-\tau)\tau^{\alpha-1}E_{\alpha,\alpha}(-\mu_n \tau^\alpha) d\tau\right| \; |\psi_n(x)| \\ \nonumber
&\leq &\frac{|g(0)|}{\Gamma(\alpha)}\sum_{n=1}^{\infty} |f_n | \;\|\psi_n\|_{C(\overline{\Omega})}\;t^{\alpha-1}+\|g\|_{C^1[0,T]} \sum_{n=1}^{\infty}\frac{1}{\mu_n}|f_n|\; \|\psi_n\|_{C(\overline{\Omega})}\\ \nonumber
&\leq&C \frac{|g(0)|}{\Gamma(\alpha)}\left( \sum_{n=1}^{\infty}\frac{1}{\mu_n^{2m}}\right)^{1/2} \; \left( \sum_{n=1}^{\infty}\mu_n^{2(k+m)}f_{n} ^2\right)^{1/2} t^{\alpha-1}+\nonumber\\
&&C \|g\|_{C^1[0,T]} \left( \sum_{n=1}^{\infty}\frac{1}{\mu_n^{2m}}\right)^{1/2} \;
\left( \sum_{n=1}^{\infty}\mu_n^{2(m+k-1)}f_{n} ^2\right)^{1/2} \nonumber\\
&\leq & C  \|f\|_{\mathcal{D}((-\Delta)^{\gamma})}t^{\alpha-1} +C \|g\|_{C^1[0,T]} \|f\|_{ \mathcal{D}((-\Delta)^{\gamma-1})}.
\end{eqnarray}

For any fixed $t_0\in (0,T]$, we can always choose $\epsilon_0>0$ such that $t_0>\epsilon_0$ and the series defining $\partial_t u_1, \partial_t u_2$ are uniformly convergent in $[\epsilon_0, T]$, and consequently the series are continuous in $[\epsilon_0, T]$, especially at $t_0>\epsilon_0$. Since $t_0>0$ is arbitrary, we know that $\partial_t u_1(x,\cdot), \partial_t u_2(x,\cdot)$ are continuous in $(0, T]$.

As for $ \partial_t u_3(x,t):= \partial_t z_b(x,t)\equiv  \partial_t z_c(x,t)+ \partial_t \Lambda_b(x,t)$ defining in (\ref{AA3}), since $ \partial_t \Lambda_b(x,\cdot)\in C[0,T]$ with  $b(x,\cdot)\in C[0,T]$ for every $x\in\Omega$, let us show $  z_c(x,\cdot)\in C^1(0,T]$. For the series defining $\partial_t z_c(x,t)$,  using the boundedness of
$0<E_{\alpha,1}(-\mu_n t^\alpha)\leq 1$, we also have the estimate
\begin{eqnarray}\label{U3t}
& &\sum_{n=1}^{\infty}|\left( \partial_t\Lambda_b(\cdot,0),\psi_n(\cdot)\right)|\;|  E_{\alpha,1}(-\mu_n t^\alpha)|  \; |\psi_n(x)|+ \nonumber\\
&&\sum_{n=1}^{\infty}  \left| \int_0^t \left( \partial_t^2\Lambda_b(\cdot,t-\tau),\psi_n(\cdot)\right)  E_{\alpha,1}(-\mu_n \tau^\alpha)  d\tau\right|\; |\psi_n(x)|
\nonumber\\&\leq&
\sum_{n=1}^{\infty}  \left|\left( \partial_t\Lambda_b(\cdot,0),\psi_n(\cdot)\right)_{L^2(\Omega)}   \right| \; \|\psi_n\|_{C(\overline{\Omega})}+
\nonumber\\& &
T^{1/2}\sum_{n=1}^{\infty} \|(\partial_t^2\Lambda_b(*,\cdot),\psi_n(*))_{L^2(\Omega)}\|_{L^2(0,T)}  \;
 \|\psi_n\|_{C(\overline{\Omega})}\nonumber\\
 &\leq& C \left( \sum_{n=1}^{\infty}\frac{1}{\mu_n^{2m}}\right)^{1/2} \; \left( \sum_{n=1}^{\infty}\mu_n^{2(k+m)}\left|(\partial_t\Lambda_b(\cdot,0),\psi_n(\cdot)\right)_{L^2(\Omega)}| ^2\right)^{1/2}+\nonumber\\
&&C \left( \sum_{n=1}^{\infty}\frac{1}{\mu_n^{2m}}\right)^{1/2} \;
\left( \int_0^T\sum_{n=1}^{\infty} \mu_n^{2(k+m)} |(\partial_t^2\Lambda_b(\cdot,\tau),\psi_n(\cdot))_{L^2(\Omega)}|^2 d\tau\right)^{1/2}\nonumber\\
&\leq& C \|\partial_t\Lambda_b\|_{C([0,T]; \mathcal{D}((-\Delta)^{\gamma}))}+C \|\partial_t^2\Lambda_b\|_{L^2(0,T; \mathcal{D}((-\Delta)^{\gamma}))}.
  \end{eqnarray}
So, for $b\in C^2([0,T]; H^{5/2}(\partial\Omega))$, the series defining $\partial_t  z_c(x,\cdot)$ for every $x\in \Omega$ are absolutely uniformly convergent in $(0,T]$ in terms of (\ref{U3t}) and (\ref{liu-3333}) for $l=4$, which ensure that $\partial_t u_3$, and consequently $\partial_t u(x,\cdot)=\partial_t u_1+\partial_t u_2+\partial_t u_3$, are continuous in  $(0,T]$ for every $x\in \Omega$, which verify $u(x,\cdot)\in C^1(0,T]$ for every $x\in \Omega$.
By taking $\gamma=2$, we also have the estimate from (\ref{U1t})-(\ref{U3t}) that
\begin{eqnarray}\label{REvis1}
\| \partial_t u(\cdot,t)\|_{C(\overline{\Omega})}
&\leq& C \|u_0\|_{\mathcal{D}((-\Delta)^{3})}t^{\alpha-1}+C \|f\|_{\mathcal{D}((-\Delta)^{2})}t^{\alpha-1} +\nonumber\\& &C \|g\|_{C^1[0,T]} \|f\|_{ \mathcal{D}((-\Delta)^{1})}
+C \|\partial_t b\|_{C([0,T];H^{5/2}(\partial\Omega))}+\nonumber\\& &C \|\partial_t^2 b\|_{L^2(0,T; H^{5/2}(\partial\Omega))},
\end{eqnarray}
which says $\partial_tu \in L^1((0,T), C(\overline\Omega))$. So we have verified $u(x,\cdot)\in W_t^1(0,T)$ for every $x\in\Omega$, i.e., $u\in S_1(\Omega_T)$.

Finally, We will prove $u(\cdot,t)\in C^2(\Omega)$ for any $t>0$. It is enough to prove that $u_i(\cdot,t), i=1,2,3$ is twice continuously differentiable with respect to  $x\in\Omega$ for any fixed $t\in[0,T]$.
For $\Omega\subset \mathbb{R}^d$ with $d=1,2,3$, we always have $H^4(\Omega)\subset C^2(\overline\Omega)$. Therefore
For the series defining the 2-times derivative  $\partial_{x_i x_j} u_1$ and $\partial_{x_i x_j} u_2$, we have the estimates with $k=2$ in (\ref{AA00-3}) that
\begin{eqnarray}\label{U1-C2}
& &\sum_{n=1}^{\infty} |u_{0n} \; E_{\alpha,1}(-\mu_n t^\alpha)| \; |\partial_{x_i x_j}\psi_n(x)|\nonumber\\
&\leq& \sum_{n=1}^{\infty} |u_{0n}| \;\|\psi_n\|_{C^2(\overline{\Omega})}\leq C\sum_{n=1}^{\infty} |u_{0n}| \;\|\psi_n\|_{H^4(\Omega)}\nonumber\\
&\le& C \left( \sum_{n=1}^{\infty}\frac{1}{\mu_n^{2m}}\right)^{1/2} \; \left( \sum_{n=1}^{\infty}\mu_n^{2(2+m)}u_{0n} ^2\right)^{1/2}\le C\|u_0\|_{\mathcal{D}((-\Delta)^\gamma)}
\end{eqnarray}
and
\begin{eqnarray}\label{U2-C2}
& &\sum_{n=1}^{\infty} |f_n \; g_n(t)| \; |\partial_{x_i x_j}\psi_n(x)|
\nonumber\\&\leq& \|g\|_{C[0,T]} \sum_{n=1}^{\infty}\frac{1}{\mu_n}|f_n|\; \|\psi_n\|_{C^2(\overline{\Omega})}
\leq C \|g\|_{C[0,T]} \sum_{n=1}^{\infty}\frac{1}{\mu_n}|f_n|\; \|\psi_n\|_{H^4(\Omega)}
\nonumber\\
&\leq&
C\|g\|_{C[0,T]} \left( \sum_{n=1}^{\infty}\frac{1}{\mu_n^{2m}}\right)^{1/2}\;
\left( \sum_{n=1}^{\infty}\mu_n^{2(2+m-1)}f_{n} ^2\right)^{1/2}
\nonumber\\&\leq&
C\|g\|_{C[0,T]}
\|f\|_{\mathcal{D}((-\Delta)^{\gamma-1})}.
\end{eqnarray}
For $u_3(x,t):=z_b(x,t)\equiv z_c(x,t)+\Lambda_b(x,t)$ defining in (\ref{AA3}) which solves (\ref{FB4}), since $\Lambda_b(\cdot,t)\in C^2(\Omega)$ with
$$\|\Lambda_b(\cdot, t)\|_{C^2(\overline\Omega)}\le \|\Lambda_b(\cdot,t)\|_{H^4(\Omega)}\le C\|b(\cdot,t)\|_{H^{5/2}(\partial\Omega)}$$
for $t\in [0,T]$,  let us verify $z_c(\cdot,t)\in C^2(\Omega)$. For the series defining $\partial_{x_i x_j}z_c(x,t)$, we also have
\begin{eqnarray}\label{U3-C2}
& &\sum_{n=1}^{\infty} \;  \int_0^t \left| (\partial_t\Lambda_b(\cdot,t-\tau),\psi_n(\cdot))E_{\alpha,1}(-\mu_n \tau^\alpha) \right| d\tau \; |\partial_{x_i x_j}\psi_n(x)|\nonumber\\
&\leq&\sum_{n=1}^{\infty} \|\psi_n\|_{C^2(\overline{\Omega})}\left(\int_0^t (\partial_t\Lambda_b(\cdot,t-\tau),\psi_n(\cdot))^2 d\tau\right)^{1/2}
\nonumber\\
&\leq&C\sum_{n=1}^{\infty} \|\psi_n\|_{H^4(\Omega)}\left(\int_0^t (\partial_t\Lambda_b(\cdot,\tau),\psi_n(\cdot))^2 d\tau\right)^{1/2}
\nonumber\\
&\le&
C\left( \sum_{n=1}^{\infty}\frac{1}{\mu_n^{2m}}\right)^{1/2} \;\left( \sum_{n=1}^{\infty}\mu_n^{2(2+m)}\int_0^T(\partial_\tau\Lambda_b(\cdot,\tau),\psi_n(\cdot))^2d\tau\right)^{1/2}
\nonumber\\
&\le&C\|\partial_t\Lambda_{ b}\|_{L^2(0,T;D((-\Delta)^\gamma))}.
\end{eqnarray}
Here, we take $\gamma=k+m=3$ with $k=2>d/4, m=1>d/4$ for $d=1,2,3$. So we have the estimate
\begin{eqnarray}\label{SEP20-101}
 \nonumber \|\partial_t\Lambda_b\|_{L^2(0,T; \mathcal{D}((-\Delta)^{\gamma})}^2&\equiv&
\int_0^T \sum_{n} \mu_n^{6}(\partial_t\Lambda_b(\cdot,t),\psi_n)^2dt\\
\nonumber
&\leq& C \|\partial_t\Lambda_b\|^2_{L^2(0,T;H^6(\Omega))}
\leq C \|\partial_tb\|^2_{L^2(0,T;H^{9/2}(\partial\Omega))}
\end{eqnarray}
by (\ref{liu-3333}) for $l=6$.
Hence, for $(u_0,f,g,b)$ satisfying the specified regularities, the three series defining $\partial_{x_i x_j} u_1(\cdot,t),\partial_{x_i x_j} u_2(\cdot,t),\partial_{x_i x_j} z_c(\cdot,t)$ for every $t\in (0,T]$ are absolutely uniformly convergent in $\Omega$ in terms of (\ref{U1-C2})-(\ref{U3-C2}), which ensure that $\partial_{x_i x_j} u_1$,$\partial_{x_i x_j} u_2$, $\partial_{x_i x_j} u_3$, and consequently $\partial_{x_i x_j} u(\cdot,t)=\partial_{x_i x_j} u_1+\partial_{x_i x_j} u_2+\partial_{x_i x_j} u_3$, are continuous functions in  $\Omega$ for every $t\in [0,T]$, which verify $u(\cdot,t)\in C^2(\Omega)$ for every $t\in [0,T]$.
We also have the estimate
\begin{eqnarray}\label{REvis2}
\nonumber \| \partial_{x_i x_j} u(\cdot,t)\|_{C(\Omega)}
&\leq&  C\|u_0\|_{\mathcal{D}((-\Delta)^\gamma)}
+\|g\|_{C[0,T]}\|f\|_{\mathcal{D}((-\Delta)^{\gamma-1})}+
\\ \nonumber
&&C\|\partial_t\Lambda_{ b}\|_{L^2(0,T;D((-\Delta)^\gamma))}+\|\Lambda_b(\cdot,t)\|_{C^2(\overline\Omega)}\\ \nonumber
&\leq&  C\|u_0\|_{\mathcal{D}((-\Delta)^3)}
+\|g\|_{C[0,T]}\|f\|_{\mathcal{D}((-\Delta)^{2})}\\ \nonumber
&&C \|\partial_tb\|_{L^2(0,T;H^{9/2}(\partial\Omega))}+C\|b\|_{C([0,T], H^{5/2}(\partial\Omega))}.
\end{eqnarray}
So $u\in S_2(\Omega_T)$. The proof is complete.
\end{proof}

\section{The maximum principle}

Now we study the maximum principle of the initial-boundary problem (\ref{eq1}) with the strong solution  in the sense of definition \ref{DefReGU11}. This maximum principle plays a key role for the inverse problems of recovering $\lambda(x)$, which will be studied in our further work. In other words, we assume that  $u_0\in \mathcal{D}((-\Delta)^{3})$, $g\in C^1[0,T]$, $f \in \mathcal{D}((-\Delta)^{2})$ and $b \in C^2([0,T]; H^{9/2}(\partial\Omega))$ in this section.

For $\alpha\in (0,1)$, define $\mathrm{x}:=(x,t)$ and
$$L^\alpha u:=\partial_{0+}^\alpha u -\Delta u.$$
We need to firstly establish the Hopf Lemma for the continuous function satisfying the time fractional order equation in (\ref{eq1}), which is a generalization of the corresponding result for classic parabolic system  \cite{Friedman}.

\begin{lemma}\label{lem1}
Assume $L^\alpha u =F(x,t)\leq 0(\geq 0)$ and $u\in C(\overline\Omega_T)$ attains its maximum $M$ (minimum $m$) at point $\mathrm{x}_0=(x_0,t_0)\in \partial\Omega \times (0,T]$. If $\Omega_T$ satisfies an interior strong sphere property at $\mathrm{x}_0$ and there exists a neighbourhood $V$ of $\mathrm{x}_0$ such that $u<M$($u>m$) in
$\Omega_T\cap V$. Then it holds
\begin{equation}\label{liu4-1}
\frac{\partial u}{\partial \nu}\left.\right|_{\mathrm{x}_0}>0(<0),
\end{equation}
where $\nu$ is the outward normal direction at $\mathrm{x}_0$.
\end{lemma}

\begin{proof}
We only prove the case $F(x,t) \leq 0$, the other case can be shown analogously.
Since $\Omega_T$ satisfies an interior strong sphere  property at $\mathrm{x}_0$, there exists a ball $B=B_R(\mathrm{\overline{x}})\subset \overline{\Omega}_T$ small enough specified by
$|x-\overline{x}|^2+(t-t_0)^{2}\leq R^2$
with center $\mathrm{\overline{x}}=(\overline{x},t_0)$  such that the interior of $B$ lies in $\Omega_T\cap V$.

Since $\overline x\neq x_0$, there exists a hyperplane $\pi$ which divides the $(x,t)$-space into two half-spaces $\pi^-$ and $\pi^+$ such that
$(\overline{x},t_0)\in \pi^-$, $(x_0,t_0)\in \pi^+$, $B^+:=\pi^+\cap B\neq \varnothing$, and $|x-\overline{x}|\geq C >0$ for all $(x,t)\in B^+$. The boundary of $B^+$ consists of one part $C_1$ lying on $\partial B$ and another part $C_2$ lying on $\pi$. Obviously $|x_0-\overline{x}|^2=R^2$ due to $(x_0,t_0)\in\partial B$.

Introduce an auxiliary function
\begin{eqnarray*}
 h(x,t):=E_{\alpha,1}(-\mu (t-t_0)^{2\alpha})(e^{-\mu |x-\overline{x}|^{2}}-e^{-\mu (R^{2}-(t-t_0)^{2})})
\end{eqnarray*}
with a constant parameter $\mu>0$, which meets $h(x,t)|_{C_1}=0, h(x_0,t_0)=0$ and $h>0$ on $\overline{B^+}\setminus {C_1}$ by the monotonic decreasing of $e^{-\mu t}$ and $0<E_{\alpha,1}(-\mu (t-t_0)^{2\alpha})<1$. Now we show $L^\alpha h<0$ in $\overline{B^+}$ for $\mu >>1$. In fact, we have
\begin{eqnarray}\label{B10}
h_t(x,t)&=&
-2\mu(t-t_0)^{2\alpha-1} E_{\alpha,\alpha}(-\mu (t-t_0)^{2\alpha})\left(e^{-\mu |x-\overline{x}|^{2}}-e^{-\mu (R^{2}-(t-t_0)^{2})}\right) \nonumber \\
&&-2\mu (t-t_0) E_{\alpha,1}(-\mu (t-t_0)^{2\alpha})e^{-\mu (R^{2}-(t-t_0)^{2})},
\end{eqnarray}
where $(t-t_0)^{2\alpha-1}:=(t-t_0)^{2\alpha}(t-t_0)^{-1}$ for $t\not=t_0$. Obviously, $h_t(x,t)$ is singular at $t=t_0$, with the estimate for $(x,t)\in \overline{B^+}$ that
$$|h_t(x,t)|\le 2\mu(2+T^{2(1-\alpha)}) |(t-t_0)^{2\alpha-1}|e^{-\mu |x-\overline x|^2}\tilde E_{\alpha}(-\mu (t-t_0)^{2\alpha})$$
by $R^{2}-(t-t_0)^{2}\ge |x-\overline x|^2$, where
$$\tilde E_{\alpha}(-\mu (t-t_0)^{2\alpha}):=E_{\alpha,\alpha}(-\mu (t-t_0)^{2\alpha})+E_{\alpha,1}(-\mu (t-t_0)^{2\alpha}).$$
Due to this representation, we can verify
$\partial_{0+}^\alpha h(x,t)$
is well-defined for $t\in (0,T]$. In fact, since $E_{\alpha,\alpha}(-\mu (\cdot-t_0)^{2\alpha}), E_{\alpha,1}(-\mu (\cdot-t_0)^{2\alpha})\in C[0,T]$,
it is enough to verify the integrability of
\begin{eqnarray}\label{star1}
H(t,t_0):=\int_0^t \frac{(\tau-t_0)^{2\alpha}}{(t-\tau)^\alpha(\tau-t_0)}d \tau
\end{eqnarray}
for $t\in (0,T]$, which can be seen from
\begin{eqnarray*}
H(t,t_0)=
\begin{cases}
-\int_0^{t} \frac{(t_0-\tau)^{2\alpha}}{(t-\tau)^\alpha(t_0-\tau)}d \tau, &0<t<t_0,\\
-\int_0^{t_0} \frac{1}{(t_0-\tau)^{1-\alpha}}d \tau, &t=t_0,\\
-\int_0^{t_0} \frac{(t_0-\tau)^{\alpha}}{(t-\tau)^\alpha(t_0-\tau)^{1-\alpha}}d \tau
+\int_{t_0}^t \frac{(\tau-t_0)^{\alpha}}{(t-\tau)^\alpha(\tau-t_0)^{1-\alpha}}d \tau, &t_0<t\le T
\end{cases}
\end{eqnarray*}
clearly. Moreover, we have for $(x,t)\in \overline {B^+}$ that
\begin{eqnarray*}\label{partial1}
& &|\partial_{0+}^\alpha h(x,t)|
\nonumber\\&\le&  \frac{2\mu (2+T^{2(1-\alpha)}) e^{-\mu |x-\overline x|^2}}{\Gamma(1-\alpha)}
\int_0^t \frac{(\tau-t_0)^{2\alpha}}{(t-\tau)^\alpha|\tau-t_0|}\tilde E_{\alpha}(-\mu (\tau-t_0)^{2\alpha})d\tau
\end{eqnarray*}
by (\ref{B10}). Now we estimate the integral in the right-hand side for two cases.

Case 1: $0<t\le t_0$. Then there exists $0<\xi_t<t\le t_0$ such that
\begin{eqnarray*}
& &\int_0^t \frac{(\tau-t_0)^{2\alpha}}{(t-\tau)^\alpha|\tau-t_0|}\tilde E_{\alpha}(-\mu (\tau-t_0)^{2\alpha})d\tau
\nonumber\\&=&\tilde E_{\alpha}(-\mu (\xi_t-t_0)^{2\alpha})
\int_0^t \frac{(\tau-t_0)^{2\alpha}}{(t-\tau)^\alpha|\tau-t_0|}d\tau
\nonumber\\&\le&
C(T,t_0,\alpha)\tilde E_{\alpha}(-\mu (\xi_t-t_0)^{2\alpha})
\end{eqnarray*}
by (\ref{star1}).

Case 2: $0<t_0<t\le T$. It follows analogously that
\begin{eqnarray*}
& &\int_0^t \frac{(\tau-t_0)^{2\alpha}}{(t-\tau)^\alpha|\tau-t_0|}\tilde E_{\alpha}(-\mu (\tau-t_0)^{2\alpha})d\tau
\nonumber\\&=&
\left[\int_0^{t_0}+\int_{t_0}^t\right] \frac{(\tau-t_0)^{2\alpha}}{(t-\tau)^\alpha|\tau-t_0|}\tilde E_{\alpha}(-\mu (\tau-t_0)^{2\alpha})d\tau
\nonumber\\&=&
\tilde E_{\alpha}(-\mu (\xi_t^1-t_0)^{2\alpha})
\int_0^{t_0} \frac{(\tau-t_0)^{2\alpha}}{(t-\tau)^\alpha|\tau-t_0|}d\tau+
\nonumber\\& &
\tilde E_{\alpha}(-\mu (\xi_t^2-t_0)^{2\alpha})
\int_{t_0}^t \frac{(\tau-t_0)^{2\alpha}}{(t-\tau)^\alpha|\tau-t_0|}d\tau
\nonumber\\&\le&
C(T,t_0,\alpha)\sum_{i=1}^2 \tilde E_{\alpha}(-\mu (\xi_t^i-t_0)^{2\alpha})
\end{eqnarray*}
with $t_0\not=\xi_t^i<t$ for $i=1,2$ by (\ref{star1}). Combining these two cases together, we know for $t\in (0,T]$ that
\begin{eqnarray}\label{liu5}
\int_0^t \frac{(\tau-t_0)^{2\alpha}}{(t-\tau)^\alpha|\tau-t_0|}\tilde E_{\alpha}(-\mu (\tau-t_0)^{2\alpha})d\tau
\simeq
\tilde E_{\alpha}(-\mu (\tilde \xi_t-t_0)^{2\alpha}), \quad t\in (0,T)
\end{eqnarray}
for some $ t_0 \not =\tilde\xi_t \in (0,T)$ when $\mu\to\infty$. On the other hand, we have
\begin{eqnarray}\label{BZM1}
\Delta h=h_{x_1x_1}+h_{x_2x_2}
=
E_{\alpha,1}(-\mu (t-t_0)^{2\alpha}) e^{-\mu |x-\overline{x}|^{2}}(4\mu^2|x-\overline{x}|^2-4\mu)
\end{eqnarray}
from straightforward computations. So we have for $(x,t)\in B^{+}$ from
$$L^\alpha h\le
|\partial_{0+}^\alpha h(x,t)|-e^{-\mu |x-\overline{x}|^{2}}E_{\alpha,1}(-\mu (t-t_0)^{2\alpha})
(4\mu^2|x-\overline{x}|^2-4\mu)$$
and (\ref{partial1})-(\ref{BZM1}) that
\begin{eqnarray}\label{BZM2}
& &e^{\mu |x-\overline x|^2}\mu^{-1}L^\alpha h(x,t)
\nonumber\\&\le& C(\alpha,T)
\int_0^t \frac{(\tau-t_0)^{2\alpha}}{(t-\tau)^\alpha|\tau-t_0|}\tilde E_{\alpha}(-\mu (\tau-t_0)^{2\alpha})d\tau-
\nonumber\\& &
4E_{\alpha,1}(-\mu (t-t_0)^{2\alpha})
(\mu|x-\overline{x}|^2-1)
\nonumber\\&\le &
C(\alpha,T,t_0)\tilde E_{\alpha}(-\mu (\tilde \xi_t-t_0)^{2\alpha})-
4E_{\alpha,1}(-\mu (t-t_0)^{2\alpha})(\mu|x-\overline{x}|^2-1)
\nonumber\\&<&0
\end{eqnarray}
for $\mu>>1$ due to $|x-\overline x|_{x\in \overline{B^+}}\ge C>0$, noticing that
\begin{eqnarray*}
\tilde E_{\alpha}(-\mu (\tilde \xi_t-t_0)^{2\alpha})=
\begin{cases}
O(E_{\alpha,1}(-\mu (t-t_0)^{2\alpha})), &t\not=t_0,\\
o(E_{\alpha,1}(-\mu (t-t_0)^{2\alpha})), &t=t_0
\end{cases}
\end{eqnarray*}
as $\mu\to\infty$,
i.e., $L^\alpha h<0$ on $\overline{B^+}$ for large $\mu>0$.

For $\epsilon>0$ sufficiently small and the above function $h(x,t)$, consider  the function $v:=u+\epsilon h$. Since
$$v|_{C_2}<M, \quad v|_{C_1}=u|_{C_1}+\epsilon h|_{C_1}= u|_{C_1}\le M$$
due to $u|_{C_2}<M, h|_{C_1}\equiv 0$ and $v(\mathrm{x}_0)=u(\mathrm{x}_0)=M$, $v$ takes its maximum value $M$ in $\partial B^+$ at point $\mathrm{x}_0=(x_0,t_0)\in C_1$. On the other hand, we have from (\ref{BZM2}) that
$$L^\alpha v=L^\alpha u+\epsilon L^\alpha h=F(x,t)+\epsilon L^\alpha h\le \epsilon L^\alpha h<0 \hbox{ in } B^+,$$
which means that  $v$ cannot attain its  maximum in $\overline{B^+}$ in its interior points by \cite{Roscani} (Prop. 1), i.e., $v<M$ in the interior of $B^+$, which yields
$\frac{\partial v}{\partial \nu}|_{\mathrm{x}_0} \geq 0$.
Noticing
\begin{eqnarray*}
\frac{\partial h}{\partial \nu}|_{\mathrm{x}_0} =-2\mu e^{-\mu R^{2}} \sum_{i=1}^{2}({x_0}_i-\overline{x}_i)\frac{{x_0}_i-\overline{x}_i}{R}<0,
\end{eqnarray*}
we finally have
\begin{eqnarray*}
\frac{\partial u}{\partial \nu}|_{\mathrm{x}_0}= \frac{\partial v}{\partial \nu}|_{\mathrm{x}_0} -\epsilon \frac{\partial h}{\partial \nu}|_{\mathrm{x}_0}>0.
\end{eqnarray*}
The proof is complete.
\end{proof}

The direct application of this result leads to the following weak maximum principle for our slow diffusion system (\ref{eq1}) with impedance boundary condition, based on Theorem 3 in \cite{Luchko}.

\begin{theorem}\label{lem00}
Let $u_0\geq 0 $, $F(x,t)\geq 0$, $ b(x,t)\geq 0$, $0<\lambda^-<\lambda(x)\leq \lambda^+$ and $u(x,t)$ be the classical solution of (\ref{eq1}) defined above. Then $u(x,t)\geq 0$ in $\overline{\Omega}_T:=\overline{\Omega}\times[0,T]$.
\end{theorem}
\begin{proof}
Since $F(x,t)\geq0$, then either $u\geq0$ in $\overline{\Omega}_T$ or $u$ attains its negative minimum on
$\partial \Omega \times (0,T] \cup \overline{\Omega} \times\{0\}$ by \cite{Luchko}.  If $u\geq 0$ in $\overline{\Omega}_T$  is false, then there must exist a point $(x^*,t^*)\in \partial\Omega\times(0,T]$ such that $u(x^*,t^*)=\min_{\overline{\Omega}_T}u[\lambda,f](x,t)<0$ due to $u(x,0)=u_0\geq 0$. Then the Hopf Lemma \ref{lem1} says $\partial_\nu u(x,t)|_{(x^*,t^*)}<0$. So we have
$$\partial_\nu u(x^*,t^*)+\lambda(x)u(x^*,t^*)<0$$
from the impedance boundary condition,
which contradicts the condition
$$\partial_\nu u(x,t)+\lambda(x)u(x,t)=b(x,t)\geq 0, ~~ (x,t)\in\partial\Omega\times (0,T]. $$
Therefore, $u\geq0$ in $\overline{\Omega}_T$.
The proof is complete.
\end{proof}

Based on this weak maximum principle, there holds the following strong maximum principle.
\begin{theorem}\label{lem0}
Let $u_0> 0 $, $F(x,t)\geq 0$, $ b(x,t)\geq 0$, $0<\lambda^-<\lambda(x)\leq \lambda^+$ and $u(x,t)$ be the classical solution of (\ref{eq1}). Then, $u(x,t)>0$ for $(x,t)\in  \overline{\Omega}_T$.
\end{theorem}

\begin{proof}
Firstly, we will prove $u >0$ in $[0,T] \times \Omega$.  By the linear superposition principle, we decompose $u(x,t)=u_1(x,t)+u_2(x,t)+u_3(x,t)$ with
\begin{eqnarray}\label{AB1}
\begin{cases}
\partial_{0+}^\alpha u_1(x,t)-\Delta u_1(x,t)=0, &(x,t)\in \Omega_T,\\
\frac{\partial u_1(x,t)}{\partial \nu}+\lambda(x)u_1(x,t)=0, &(x,t)\in \partial\Omega\times [0,T],\\
u_1(x,t)(x,0)=u_0(x), &x\in\Omega,
\end{cases}
\end{eqnarray}
\begin{eqnarray}\label{AB2}
\begin{cases}
\partial_{0+}^\alpha u_2(x,t)-\Delta u_2(x,t)=F(x,t), &(x,t)\in \Omega_T,\\
\frac{\partial u_2(x,t)}{\partial \nu}+\lambda(x)u_2(x,t)=0, &(x,t)\in \partial\Omega\times [0,T],\\
u_2(x,t)(x,0)=0, &x\in\Omega,
\end{cases}
\end{eqnarray}
and
\begin{eqnarray}\label{AB1}
\begin{cases}
\partial_{0+}^\alpha u_3(x,t)-\Delta u_3(x,t)=0, &(x,t)\in \Omega_T,\\
\frac{\partial u_3(x,t)}{\partial \nu}+\lambda(x)u_3(x,t)=b(x,t), &(x,t)\in \partial\Omega\times [0,T],\\
u_3(x,t)(x,0)=0, &x\in\Omega.
\end{cases}
\end{eqnarray}
From the above weak maximum principle, we can see that $u_i(x,t)\geq 0$, $i=1,2,3$ in
$\overline{\Omega}_T$. Now we prove  $u_1>0$ in $\overline\Omega_T$ which leads to $u=u_1+u_2+u_3>0$ in $\overline\Omega_T$.

Obviously $u_1(x,0)=u_0(x)>0$. If $u_1(x^*, t^*)=0$ for some $(x^*, t^*)\in \partial\Omega\times (0,T]$, then $u_1$ takes its minimum value on the boundary point $(x^*, t^*)$. So by Lemma 4.1, it follows that $\partial_\nu u_1(x^*,t^*)<0$, leading to
$$\partial_\nu u_1(x^*,t^*)+\lambda(x^*)u_1(x^*,t^*)=\partial_\nu u_1(x^*,t^*)<0,$$
contradicting the boundary condition for $u_1(x,t)$.

If $u_1(x^*, t^*)=0$ for some point $(x^*, t^*)\in\Omega\times (0,T]$, then $u_1$ takes its minimum value $0$ in interior point. Then we have $\Delta u_1(x^*,t^*)>0$ and
$\partial_t^\alpha u_1(x^*,t^*)\le 0$ by Theorem 1 in \cite{Luchko}. So we have
$$\partial_t^\alpha u_1(x^*,t^*)-\Delta u_1(x^*,t^*)<0,$$
contradicting the equation for $u_1(x,t)$.
The proof is complete.
\end{proof}

\section{Conclusions}

In this paper, the well-posedness of the initial-boundary value problem for the time-fractional diffusion equation with inhomogenous Robin boundary condition is shown. Also, both the weak and the strong maximum principle are presented.
To the authors' knowledge, there are no literatures to identify Robin coefficient and other unknown parameters simultaneously for time-fractional diffusion equation. In the future, we will consider the inverse problem of both determining the boundary heat exchange coefficient $\lambda(x)$ and space-dependent source term $f(x)$ from the final measurement data. Such inverse problems need the maximum principles established in this work for the uniqueness of the solution.

\vskip 0.5cm

{\bf Acknowledgement:} This work is supported by  NSFC (No.11971104) and Postgraduate Research \& Practice Innovation Program of Jiangsu Province (No.SJKY 19\_0056).

\section{Appendix: Proofs}
Let us define $A: D(A) \subset L^{2}(\Omega) \rightarrow L^{2}(\Omega)$ by
$$
A u:=-\Delta u \quad \text { with } \quad D(A)=\left\{u \in H^{2}(\Omega) \mid \frac{\partial u}{\partial \nu}+\lambda(x) u=0\right\}
$$
In the sequel we will show that $A: D(A) \subset L^{2}(\Omega) \rightarrow L^{2}(\Omega)$ is a self-adjoint operator with continuous inverse and the inclusion $D(A) \subset L^{2}(\Omega)$ is compact. Then, we can conclude from (\cite{Dautray}, Theorem 6) that there exists a complete orthonormal basis of $L^{2}(\Omega)$ composed of eigenvectors of $A: D(A) \subset L^{2}(\Omega) \rightarrow L^{2}(\Omega)$
By integration by parts, it follows that
\begin{eqnarray}\label{Appendix}
(v, A u)_{L^{2}}-(u, A v)_{L^{2}}
&=&\int_{\Omega} u \Delta v-v \Delta u d x \nonumber\\
&=&\int_{\partial \Omega} u \frac{\partial v}{\partial n}-v \frac{\partial u}{\partial n} d S
=0 \quad u, v \in D(A)
\end{eqnarray}
where we have use the fact that $\frac{\partial u}{\partial n}+\lambda u=\frac{\partial v}{\partial n}+\lambda v=0$ on $\partial \Omega$ in the last step. The identity (\ref{Appendix}) implies that $A: D(A) \subset L^{2}(\Omega) \rightarrow L^{2}(\Omega)$ is symmetric. On the other hand, using integration by parts again, one has
$$
(u, A u)=\int_{\Omega}|D u|^{2} d x-\int_{\partial \Omega} \frac{\partial u}{\partial n} u d S=\int_{\Omega}|D u|^{2} d x+\int_{\partial \Omega} \lambda u^{2} d S \geq 0
$$
which ensures that it is monotone. From the classical theory of PDEs it follows that for any given $f \in L^{2}(\Omega),$ the following system
$$
\left\{\begin{array}{l}
-\Delta u+u=f \\
\frac{\partial u}{\partial n}+\lambda(x)u=0
\end{array}\right.
$$
admits a strong solution in $H^{2}(\Omega)$ (\cite{Grisvard},Theorem 2.4.2.6), which is equivalent to
$$
(A+I) u=f
$$
Hence, one has $R(I+A)=L^{2}(\Omega),$ which implies that $A: D(A) \subset L^{2}(\Omega) \rightarrow L^{2}(\Omega)$ is maximal monotone, and hence self-adjoint (\cite{Brezis}, Proposition 7.6.).

On the other hand, for any given $f \in L^{2}(\Omega)$, the elliptic system
$$
\left\{\begin{array}{l}
-\Delta u=f \\
\frac{\partial u}{\partial n}+\lambda(x)u=0 \\
\end{array}\right.
$$
admits a strong solution $u=A^{-1} f$ in $H^{2}(\Omega)$ and
$$
\left\|A^{-1} f\right\|_{H^{2}(\Omega)} \leq C\|f\|
$$
with $C>0$ independent of $f$ and $u$ (\cite{Grisvard}, Theorem 2.4.2.6). Therefore, $A: D(A) \subset L^{2}(\Omega) \subset$ $L^{2}(\Omega)$ has a continuous inverse. i.e., $A^{-1} \in \mathcal{L}\left(L^{2}(\Omega)\right) .$ Sine $A: D(A) \subset L^{2}(\Omega) \subset L^{2}(\Omega)$ is
a closed operator, it follows that $D(A)$ is a closed subspace of $H^{2}(\Omega)$ under graph norm. As the embedding $H^{2}(\Omega) \subset L^{2}(\Omega)$, we obtain that embedding $D(A) \subset L^{2}(\Omega)$ is also compact. This completes the proof.

\label{ref:ref}
 \end{document}